\documentclass[a4paper,UKenglish]{lipics}
 
\usepackage{microtype}
\usepackage[utf8]{inputenc}
\usepackage{hyperref}\hypersetup{colorlinks=false}
\usepackage{graphicx}
\newtheorem {myexample}{Example} 

\title{Gender Relations in the XIth Mexican Logic Olympiad}
\titlerunning{Gender Relations in the XIth Mexican Logic Olympiad} 

\author[1]{José Martín Castro-Manzano}
\author[2]{Verónica Reyes-Meza}
\author[3]{César López-Pérez}
\author[4]{Karen González-Fernández}
\affil[1]{Faculty of Philosophy and Humanities, UPAEP\\
  \texttt{josemartin.castro@upaep.mx}}
\affil[2]{Faculty of Psychology, UPAEP\\
  \texttt{veronica.reyes@upaep.mx}}
\affil[4]{Chair of the Logic Olympiad, AML\\
  \texttt{rudoytecnico@gmail.com}}
\affil[4]{Faculty of Philosophy, UP\\
  \texttt{karengf@gmail.com}}
\authorrunning{J.\,M. Castro-Manzano et. al.} 

\Copyright{José Martín Castro-Manzano}


\serieslogo{logo_ttl}
\volumeinfo
  {M. Antonia {Huertas}, Jo\~ao {Marcos}, Mar\'ia {Manzano}, Sophie {Pinchinat}, \\
  Fran\c{c}ois {Schwarzentruber}}
  {5}
  {4th International Conference on Tools for Teaching Logic}
  {1}
  {1}
  {17}
\EventShortName{TTL2015}

\begin{document}

\maketitle

\begin{abstract}
Science Olympiads are academical competitions with social impact due to the fact that they allow the detection of talent and promote science and critical thinking. In this work we portray the Mexican Logic Olympiad and we describe the proportion of women and men within its XIth edition with special focus on the issue of women's underrepresentation.
\end{abstract}

\section{Introduction}
Science Olympiads are academical competitions with social impact because they allow the detection of talent, due to selection criteria; and increase interest in science, both in students and professors, for they promote the further study of science and critical thinking. The relevance of Science Olympiads, thus, is not only associated with higher order education, but with fairness and social justice.
    
In Mexico the Science Olympiads cover a wide range of disciplines: mathematics,\footnote{\href{http://www.ommenlinea.org}{Olimpiada Mexicana de Matemáticas.}} chemistry,\footnote{\href{http://depa.fquim.unam.mx/olimpiada}{Olimpiada de Química del Distrito Federal.}} physics,\footnote{\href{http://rmf.fciencias.unam.mx/~onf/Olimpiadas/index.php?option=com_content&view=section&id=5&Itemid=37}{Olimpiada Mexicana de Física.}} informatics,\footnote{\href{http://www.olimpiadadeinformatica.org.mx}{Olimpiada Mexicana de Informática.}} history\footnote{\href{http://www.amc.edu.mx/amc/index.php?option=com_content&view=article&id=87&catid=40&Itemid=125}{Olimpiada Mexicana de Historia.}} and philosophy\footnote{\href{http://ipomexico.jimdo.com}{Olimpiada Mexicana de Filosofía.}} Olympiads are well known; however, in this divergent context where these sciences seem to be disjoint due to the lack of a common object of study, a unique methodology, or similar applications, there is an Olympiad that looks for unity within diversity: the Mexican Logic Olympiad (MLO),\footnote{\href{http://olimpiada.academiamexicanadelogica.org/}{Olimpiada Mexicana de Lógica.}} since its object of study is what the other sciences, albeit different, have in common: reasoning. 

Our main goal in this work is to portray the MLO and describe the proportion of women and men within its XIth edition with special focus on the issue of women's underrepresentation. In order to accomplish these goals we review the concept of logic (Section~\ref{sec:1}), we give details about how the MLO works (Section~\ref{sec:2}) and we discuss the data obtained during the MLO (Section~\ref{sec:3}). Finally, we close with a summary and remarks about future work (Section~\ref{sec:4}). 

\section{General aspects of logic}
\label{sec:1}
Reasoning is a process that produces new information given previous data by following certain norms that can be defined mathematically. These norms allow us to describe inference as the unit of measurement of reasoning: inference may be correct or incorrect depending on those norms. The science that studies such norms is Logic.

The understanding of these norms depends on three equivalent approaches: the semantical, the syntactical, and the abstract one. From the semantical standpoint the central concept of such norms is that of \textit{interpretation} and defines our notions of \textit{satisfaction}, \textit{model} and \textit{logical truth} as denoted by $\models$ and introduced by Tarski's \textit{On the Concept of Logical Consequence}~\cite{TARSKIb}. From the syntactical point of view the main concept of those norms is that of \textit{deducibility} and characterizes our intuitions of \textit{proof}, \textit{demonstration} and \textit{theorem} as denoted by $\vdash$ and introduced by Carnap's \textit{The Logical Syntax of Language}~\cite{CARNAP}. Finally, from the abstract standpoint the idea behind those norms is the \textit{relation of consequence} that generalizes the previous accounts as in Tarski's \textit{On some Fundamental Concepts of Metamathematics}, LSM, art. III~\cite{TARSKIa}, hence rendering logical consequence as the \textit{proprium} of logic.

This relation has several uses and representations. The typical use is in argumentation, the typical representation is in the argument, usually defined as a finite sequence of propositions $\phi_1,\ldots,\phi_n$ ordered in such a way that $\phi_1,\ldots,\phi_{n-1}$ are \textit{premises} and $\phi_n$, the \textit{conclusion}, is obtained by such relation of consequence. Consider a couple of examples:

\begin{myexample}\label{ex:1}
Either Jones owns a Ford, or Brown is in Boston. But Brown is not in Boston. Thus, Jones owns a Ford.
\end{myexample}

\begin{myexample}\label{ex:2}
All persons against unfair laws are rational. Every rational being is free. Thus, every free person is against unfair laws.
\end{myexample}

Example~\ref{ex:1} encodes a correct argument, Example~\ref{ex:2} does not. Justifying the (in)correction of these examples is beyond the scope of this work, but we can briefly show why Example~\ref{ex:2} is incorrect by using (standard) Logic. Let us suppose Example~\ref{ex:2} is a correct argument, then it must be impossible to find a counter-example, that is to say, it must be impossible to find a set of true premises that satisfies the argument's structure but yields a false conclusion; however, if we use the substitution set $\Theta = \lbrace$\textit{persons against  unfair laws/positive number, rational being/number greater than 0, free person/real number}$\rbrace$ we obtain true premises and a false conclusion, hence showing the argument is not correct. 

\section{Details of the MLO}
\label{sec:2}
In an effort to promote the further study of the norms of Logic in order to develop scientific and critical thinking, in 2004 a group of graduate students from the First Certificate in Logic (\textit{Primer Diplomado en Lógica}) offered by the National Autonomous University's Institute for Philosophical Research (\textit{Instituto de Investigaciones Filosóficas}) and the Mexican Academy of Sciences (\textit{Academia Mexicana de Ciencias}) (Cristian Diego Alcocer, María del Carmen Cadena Roa y Ricardo Madrid†), and thanks to the support by the Mexican Academy of Logic (\textit{Academia Mexicana de Lógica}), called the First National Logic Olympiad (\textit{Primera Olimpiada Nacional de Lógica}). Since then there have been eleven consecutive editions that have gathered students, professors, and researchers from all over the country, Argentina, Colombia, and Paraguay. Up next we give details about how the MLO works. 

\subsection{Divisions and stages}
The MLO is divided into three divisions: 

\begin{itemize}
\item High school (\textit{Bachillerato}).
\item Undergraduate (\textit{Licenciatura}).
\item Masters.
\end{itemize}

The High School and Undergraduate divisions are self-explanatory; the Masters division, introduced during the XIth edition for the first time, includes advanced undergraduate students that have been winners in previous editions and graduate students.

These divisions enter the competition in two stages:

\begin{itemize}
\item Knock-out round. 
\item Final round.
\end{itemize}

\subsection{Exams}
The exams are composed by 30 multiple choice questions designed and reviewed by the Academic Committe of the MLO composed by researchers, professors, and winners of previous editions. Since the MLO is a competition designed to measure logical skills, and not knowledge about logic, it requires logical training in three big areas that we are going to illustrate with references and toy examples (the examples, of course, do not fully express all the complexity and quality of the exams, but we hope they will give a general impression of how regular exams look like):  

\begin{itemize}
\item Classical First Order Logic.
\item Informal Logic.
\item Problem Solving.
\end{itemize}

Classical First Order Logic captures deductive reasoning skills and requires, for instance, contents similar to Copi's \textit{Symbolic Logic}~\cite{COPI}, Enderton's \textit{Mathematical Introduction to Logic}~\cite{ENDERTON}, Gamut's first volume of \textit{Logic, Language, and Meaning}~\cite{GAMUT}, Smullyan's \textit{First-Order Logic}~\cite{SMULLYANb}, or Mendelson's \textit{Introduction to Mathematical Logic}~\cite{MENDELSON}. The next one is a typical decision example:

\begin{myexample}
\textit{Since some animals are vertebrate and some cats are vertebrate, it follows that some animals are cats}. How is the previous argument?
\end{myexample}
\begin{itemize}
\item a) Valid, with true premises, and true conclusion
\item b) Invalid, with true premises, and true conclusion
\item c) Invalid, with false premises, and false conclusion
\item d) Valid, with true premises, and false conclusion
\item e) Invalid, with true premises, and false conclusion
\end{itemize}

Informal Logic tries to capture non-deductive skills related to argumentation and supposes, for instance, Perelman \& Olbrechts-Tyteca's \textit{Treatise on Argumentation}~\cite{PERELMAN}, Toulmin's \textit{The Uses of Argument}~\cite{TOULMIN} or Walton's \textit{A Pragmatic Theory of Fallacy}~\cite{WALTON}. The next one is a typical fallacy detection example:

\begin{myexample}
What fallacy best suits the next argument? \textit{If Europe is a Jupiter's satellite, then Europe spins around Jupiter. Since there are goats living in Europe, it must be the case that there are goats spinning around Jupiter}.
\end{myexample}
\begin{itemize}
\item a) False cause
\item b) Ambiguity
\item c) \textit{Ad populum}
\item d) Hasty generalization
\item e) \textit{Ad hominem}
\end{itemize}

Problem Solving tries to find combined reasoning skills to solve problems, riddles, and puzzles: it assumes contents similar to that of  Smullyan's classics such as \textit{What is the name of this book?}~\cite{SMULLYANb}; the next example, however, is due to Gardner:	

\begin{myexample}
You've got three boxes. Each one contains two balls. One box contains two black balls; another box contains two white balls; the last box contains one black ball and one white ball. The boxes are labelled according to the color of the balls it contains: BB, WW, BW. Now, someone has changed the labels of the boxes in such a way that every box is now incorrectly labelled. You can extract one ball at a time from any box without looking inside and, by this process, you must determine the content of each box. What is the smallest number of drawings needed to accomplish this task? 
\end{myexample}
\begin{itemize}
\item a) 4
\item b) 1
\item c) 3
\item d) 5
\item e) 2
\end{itemize}

\section{Results}
\label{sec:3}
The XIth edition of the MLO took place in Puebla City on may 31, 2014, at the facilities of UPAEP, and was organized by the Mexican Academy of Logic (AML) and the Faculty of Philosophy and Humanities at UPAEP. With the available data, thanks to the AML, we have pointed out some interesting facts regarding the participation of women.

\begin{figure}[h!]
  \includegraphics[width=9.5cm]{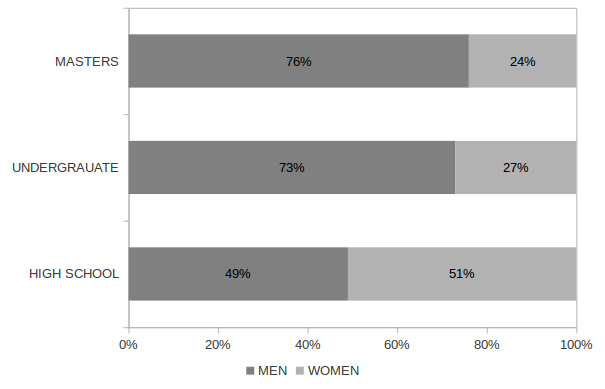}
  \caption{Proportion of gender by division}  
  \label{fig:1}
\end{figure}

The first fact we notice is that \textit{division} and \textit{gender} are inversely proportional variables: in High School the slight majority was composed by women (51\%), while in the Undergraduate and Masters division the presence was smaller, 27\%, and 24\% respectively (Figure~\ref{fig:1}).

\begin{figure}[h!]
  \includegraphics[width=9.5cm]{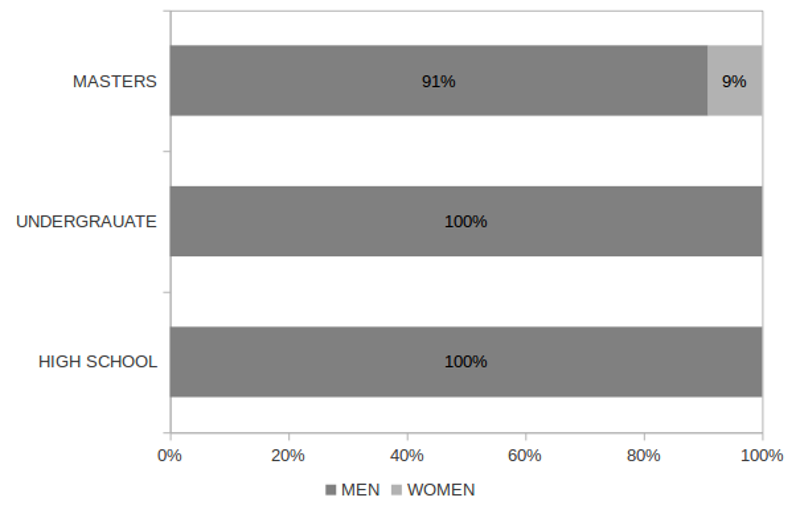}
  \caption{Proportion of medallists by division}  
  \label{fig:2}
\end{figure}

The second fact we observe is that when we consider the set of medallists (i.e, three first places) the High School and Undergraduate divisions were dominated by men (100\%), while in the Masters division a woman achieved a third place (Figure~\ref{fig:2}). Note that, since the MLO grants a medal to each and every winner, the number of medallists is greater than three; the same happens with the first ten places by division.   

\begin{figure}[h!]
  \includegraphics[width=9.5cm]{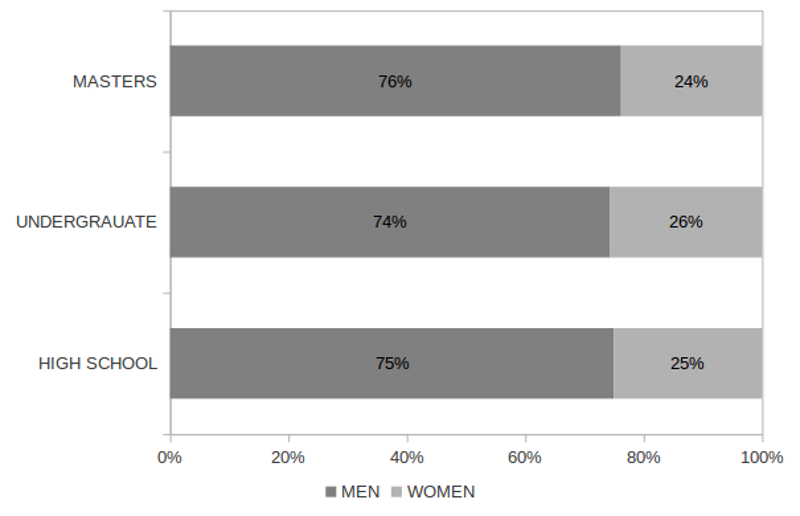}
  \caption{Proportion of first ten places by division}  
  \label{fig:3}
\end{figure}

Thirdly, given that the MLO also acknowledges the first ten places by division, we report that 25\%, 26\%, and 24\% were women, respectively (Figure~\ref{fig:3}). 

\begin{figure}[h!]
  \includegraphics[width=9.5cm]{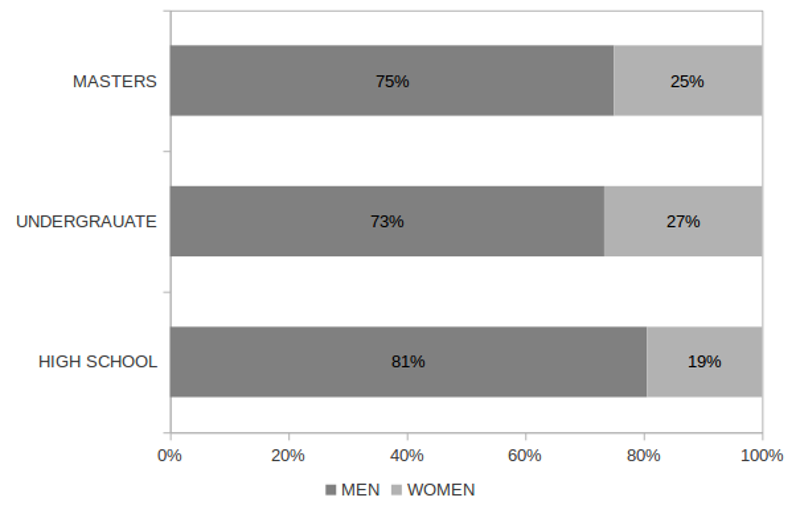}
  \caption{Proportion of coaches by division}  
  \label{fig:4}
\end{figure}

Finally since, as we said previously, Science Olympiads also target professors, we report the participation of coaches by gender and division: 19\%, 27\%, and 25\% were women, respectively (Figure~\ref{fig:4}). 

\begin{figure}[h!]
  \includegraphics[width=11cm]{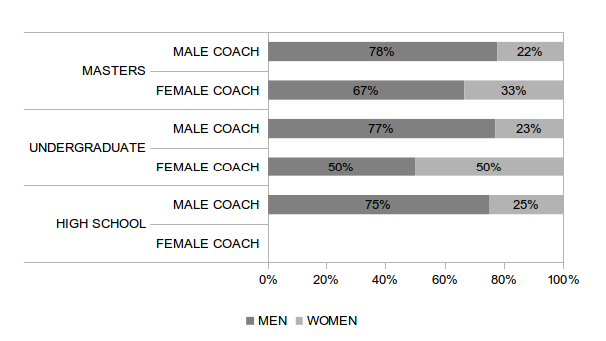}
  \caption{Proportion of coaches' gender w.r.t students' gender}  
  \label{fig:6}
\end{figure}

\begin{figure}[h!]
  \includegraphics[width=14.5cm]{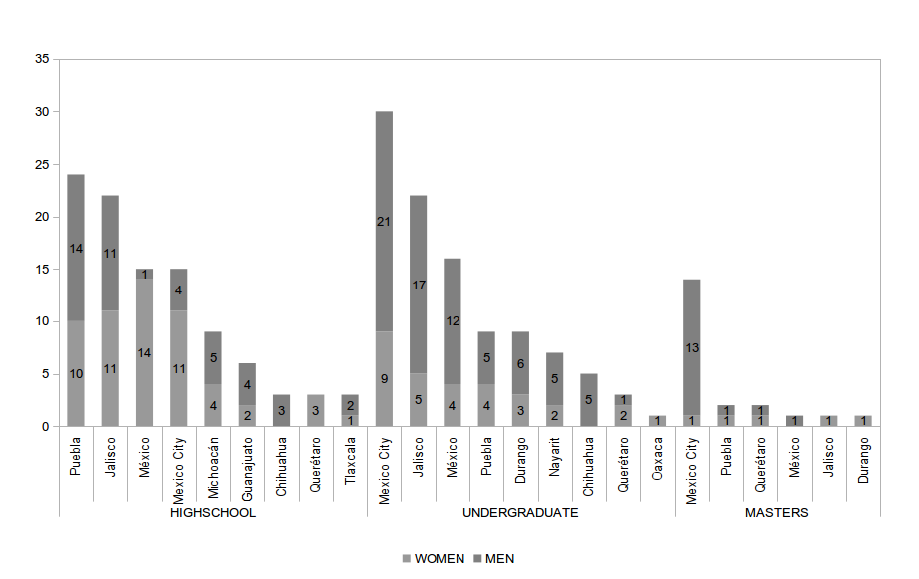}
  \caption{Proportion by state}  
  \label{fig:7}
\end{figure}

\begin{figure}[h!]
  \includegraphics[width=12cm]{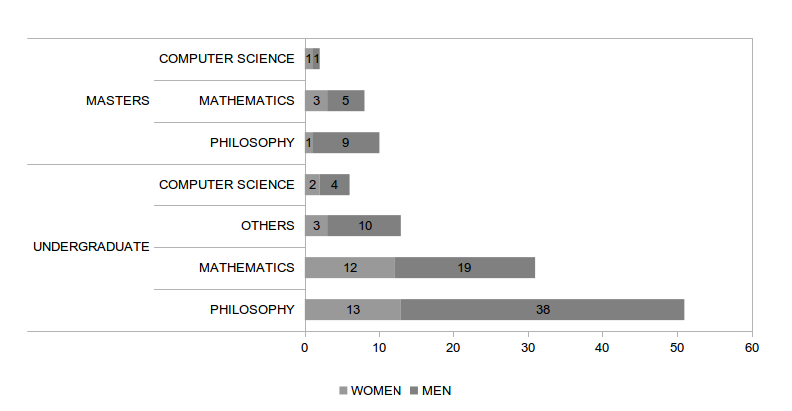}
  \caption{Proportion by major/background}  
  \label{fig:8}
\end{figure}

These descriptive results around the variables of \textit{division} and \textit{gender} are quite interesting, but we have also looked into some other variables such as the coaches' gender w.r.t students' gender (Figure~\ref{fig:6}), state (i.e., \textit{city}) (Figure~\ref{fig:7}), and current major/background (Undergraduate/Masters) (Figure~\ref{fig:8}). 

\section{Conclusion}
\label{sec:4}
Although we are currently studying more data and testing some other correlations with special focus on the participation of women, we can momentarily conclude that our report indicates that the participation of women seems to be inversely proportional to the division, despite they have been successful in the competition, both as students and coaches. Nevertheless, the cause of this underrepresentation still remains unanswered.

\begin{figure}[h!]
  \includegraphics[width=12cm]{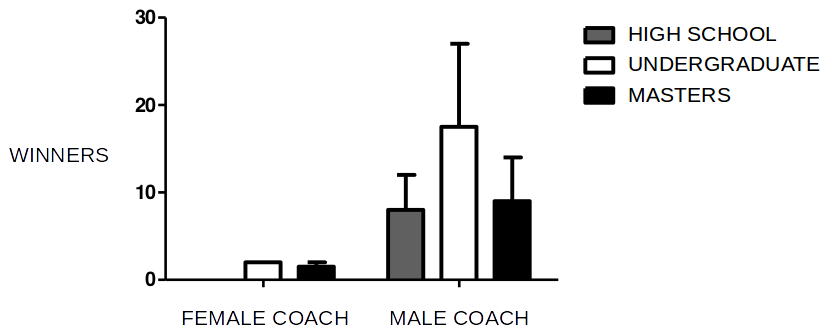}
  \caption{We have checked whether \textit{coaches' gender} has an effect in the results and it turns out this factor has significance. We assessed \textit{coaches' gender} \textit{vs division} and we observed \textit{coaches' gender} accounts for 45.98\% of the total variance giving $F = 7.31$, $DFn=1$ and $DFd=6$ with a $P\ value = 0.0354$}  
  \label{fig:9}
\end{figure}

We know there are certain disciplines in which women's participation is significantly lesser than men's (the inverse seems to be true as well), which had led some researchers to consider sexual dimorphism~\cite{MACCOBY,FRITH}, hormonal differences~\cite{HAMPSON,POSTMA,HALPERN} or environmental conditions~\cite{AMUNTS,FISHER} to provide an answer: their results, albeit interesting, are prone to objections and are quite controversial to say the less~\cite{SPELKE}. But, as our report shows, it seems clear that there is a serious loss of women when they go from High School to Undergraduate level, and we need an explanation of this phenomenon. We can also observe this loss has an evident effect on the whole competition, as our two-way ANOVA analysis shows (Figure~\ref{fig:9}), and so, it is clear that we need new policies to improve this situation.

Finally, we would like to add that the answers to Examples 3, 4, and 6 are the options \textit{b}. 

\subparagraph*{Acknowledgements}
The authors would like to thank the reviewers for their precise corrections and useful comments. Financial support given by UPAEP Grants 30108-1030 and 30108-1008.








\newpage
\thispagestyle{empty}
{\ }

\end{document}